\documentclass[twoside,11pt]{amsart}
\usepackage{amsmath,amssymb,amsthm}

\begin{document}

\setlength{\textwidth}{145mm} \setlength{\textheight}{203mm}
\setlength{\parindent}{0mm} \setlength{\parskip}{2pt plus 2pt}

\frenchspacing



\numberwithin{equation}{section}
\newtheorem{thm}{Theorem}[section]
\newtheorem{lem}[thm]{Lemma}
\newtheorem{prop}[thm]{Proposition}
\newtheorem{cor}[thm]{Corollary}
\newtheorem{probl}[thm]{Problem}

\newtheorem{defn}{Definition}[section]
\newtheorem{rem}{Remark}[section]
\newtheorem{exa}{Example}

\newcommand{\be}[1]{\begin{equation}\label{#1}}
\newcommand{\ee}{\end{equation}}


\newcommand{\X}{\mathfrak{X}}
\newcommand{\B}{\mathcal{B}}
\newcommand{\s}{\mathfrak{S}}
\newcommand{\g}{\mathfrak{g}}
\newcommand{\W}{\mathcal{W}}
\newcommand{\Lgr}{\mathrm{L}}
\newcommand{\dd}{\mathrm{d}}
\newcommand{\n}{\nabla}
\newcommand{\pd}{\partial}
\newcommand{\LL}{\mathcal{L}}

\newcommand{\diag}{\mathrm{diag}}
\newcommand{\End}{\mathrm{End}}
\newcommand{\im}{\mathrm{Im}}
\newcommand{\id}{\mathrm{id}}

\newcommand{\ie}{i.e. }
\newfont{\w}{msbm9 scaled\magstep1}
\def\R{\mbox{\w R}}
\newcommand{\norm}[1]{\left\Vert#1\right\Vert ^2}
\newcommand{\nnorm}[1]{\left\Vert#1\right\Vert ^{*2}}
\newcommand{\nN}{\norm{N}}
\newcommand{\nP}{\norm{\nabla P}}
\newcommand{\nnP}{\nnorm{\nabla P}}
\newcommand{\tr}{{\rm tr}}

\newcommand{\thmref}[1]{Theorem~\ref{#1}}
\newcommand{\propref}[1]{Proposition~\ref{#1}}
\newcommand{\corref}[1]{Corollary~\ref{#1}}
\newcommand{\secref}[1]{\S\ref{#1}}
\newcommand{\lemref}[1]{Lemma~\ref{#1}}
\newcommand{\dfnref}[1]{Definition~\ref{#1}}

\frenchspacing


\title{$P$-connection on Riemannian almost product manifolds}

\author{Dimitar Mekerov}

\maketitle

{\small
\textbf{Abstract} \\
In the present work\footnote{Partially supported by project
RS09-FMI-003 of the Scientific Research Fund, Paisii Hilendarski
University of Plovdiv, Bulgaria}, we introduce a linear connection
(preserving the almost product structure and the Riemannian
metric) on Riemannian almost product manifolds. This connection,
called $P$-connection, is an analogue of the first canonical
connection of Lichnerowicz in the Hermitian geometry and the
$B$-connection in the geometry of the almost complex manifolds
with Norden metric. Particularly, we consider the $P$-connection
on a the class of manifolds with nonintegrable almost product
structure.
}


\section{Introduction}

In \cite{Gan-Gri-Mih2}  a linear connection, called
$B$-connection, is introduced on almost complex manifolds with
Norden (or anti-Hermitian) metric. This connection (preserving the
almost complex structure and the Norden metric) is an analogue of
the first canonical connection of Lichnerowicz in the Hermitian
geometry (\cite{Gau}, \cite{Li}, \cite{Ya}). In \cite{Mek3} the
$B$-connection is considered on a class of almost complex
manifolds with Norden metric and nonintegrable almost complex
structure. This is the class $\W_3$ of the quasi-K\"ahler
manifolds with Norden metric.

In the present work, we introduce a linear connection (preserving
the almost product structure and the Riemannian metric) on
Riemannian almost product manifolds. This connection, called
$P$-connection, is an analogue of the first canonical connection
of Lichnerowicz in the Hermitian geometry and the $B$-connection
in the geometry of the almost complex manifolds with Norden
metric. Particularly, we consider the $P$-connection on the
manifolds of the class $\W_3$ from the classification in
\cite{Sta-Gri}.

The systematic development of the theory of Riemannian almost
product manifolds was started by K. Yano \cite{Ya}. In \cite{Nav}
A.~M.~Naveira gives a classification of these manifolds with
respect to the covariant differentiation of the almost product
structure. Having in mind the results in \cite{Nav}, M.~Staikova
and K.~Gribachev give in \cite{Sta-Gri} a classification of the
Riemannian almost product manifolds with zero trace of the almost
product structure.

\section{Preliminaries}

Let $(M,P,g)$ be a \emph{Riemannian almost product manifold},
\ie{} a differentiable manifold $M$ with a tensor field $P$ of
type $(1,1)$ and a Riemannian metric $g$ such that
\begin{equation}\label{Pg}
    P^2x=x,\quad g(Px,Py)=g(x,y)
\end{equation}
for arbitrary $x$, $y$ of the algebra $\X(M)$ of the smooth vector
fields on $M$. Obviously $g(Px,y)=g(x,Py)$.

Further $x,y,z,w$ will stand for arbitrary elements of $\X(M)$.

In this work we consider Riemannian almost product manifolds with
$\tr{P}=0$. In this case $(M,P,g)$ is an even-dimensional
manifold.

The classification in \cite{Sta-Gri} of Riemannian almost product
manifolds is made with respect to the tensor field $F$ of type
(0,3), defined by
\begin{equation}\label{2}
F(x,y,z)=g\left(\left(\nabla_x P\right)y,z\right),
\end{equation}
where $\nabla$ is the Levi-Civita connection of $g$. The tensor
$F$ has the following properties:
\begin{equation}\label{3}
\begin{array}{l}
    F(x,y,z)=F(x,z,y)=-F(x,Py,Pz),\\[6pt] F(x,y,Pz)=-F(x,Py,z).
\end{array}
\end{equation}
The basic classes of the classification in \cite{Sta-Gri} are
$\W_1$, $\W_2$ and $\W_3$. Their intersection is the class $\W_0$
of the \emph{Riemannian $P$-manifolds}, determined by the
condition $F(x,y,z)=0$ or equivalently $\n P=0$. In the
classification there are include the classes $\W_1\oplus\W_2$,
$\W_1\oplus\W_3$, $\W_2\oplus\W_3$ and the class
$\W_1\oplus\W_2\oplus\W_3$ of all Riemannian almost product
manifolds.

In the present work we consider manifolds from the class $\W_3$.
This class is determined by the condition
\begin{equation}\label{sigma}
    \mathop{\s}_{x,y,z} F(x,y,z)=0,
\end{equation}
where $\mathop{\s}_{x,y,z}$ is the cyclic sum by $x, y, z$. This
is the only class of the basic classes $\W_1$, $\W_2$ and $\W_3$,
where each manifold (which is not Riemannian $P$-manifold) has a
nonintegrable almost product structure $P$. This means that in
$\W_3$ the Nijenhuis tensor $N$, determined by
\begin{equation*}\label{4'}
    N(x,y)=\left(\nabla_x P\right)Py-\left(\nabla_{Px} P\right)y
    +\left(\nabla_y P\right)Px-\left(\nabla_{Py} P\right)x,
\end{equation*}
is non-zero.

Further,  manifolds of the class $\W_3$ we call \emph{Riemannian
$\W_3$-manifolds}.

As it is known the curvature tensor field $R$ of a Riemannian
manifold with metric $g$ is determined by $
    R(x,y)z=\nabla_x \nabla_y z - \nabla_y \nabla_x z -
    \nabla_{[x,y]}z
$ and the corresponding tensor field of type $(0,4)$ is defined as
follows $
    R(x,y,z,w)=g(R(x,y)z,w).
$

Let $(M,P,g)$ be a Riemannian almost product manifold and
$\{e_i\}$ be a basis of the tangent space $T_pM$ at a point $p\in
M$. Let the components of the inverse matrix of $g$ with respect
to $\{e_i\}$ be $g^{ij}$. If $\rho$ and $\tau$ are the Ricci
tensor and the scalar curvature, then $\rho^*$ and $\tau^*$,
defined by $\rho^*(y,z)=g^{ij}R(e_i,y,z,Pe_j)$ and
$\tau^*=g^{ij}\rho^*(e_i,e_j)$, are called an \emph{associated
Ricci tensor} and an \emph{associated scalar curvature},
respectively. We will use also the trace
$\tau^{**}=g^{ij}g^{ks}R(e_i,e_k,Pe_s,Pe_j)$.

The \emph{square norm} of $\nabla P$ is defined by
\begin{equation}\label{5}
\nP=g^{ij}g^{ks}g\left(\left(\nabla_{e_i}P\right)e_k,\left(\nabla_{e_j}P\right)e_s\right).
\end{equation}
Obviously $\nP=0$ iff $(M,P,g)$ is a Riemannian $P$-manifold. In
\cite{Mek1} it is proved that if $(M,P,g)$ is a Riemannian
$\W_3$-manifold then
\begin{equation}\label{6}
    \nP=-2g^{ij}g^{ks}g\left(\left(\nabla_{e_i}P\right)e_k,\left(\nabla_{e_s}P\right)e_j\right)
    =2\left(\tau-\tau^{**}\right).
\end{equation}

A tensor $L$ of type (0,4) with the pro\-per\-ties%
\begin{gather}%
L(x,y,z,w)=-L(y,x,z,w)=-L(x,y,w,z),
\label{2.11'}\medskip\\[4pt]
\mathop{\s} \limits_{x,y,z} L(x,y,z,w)=0 \quad \textit{(the first
Bianchi identity)}\label{2.12}\medskip
\end{gather}  %
is called a \emph{curvature-like tensor}. Moreover, if the
curvature-like tensor $L$ has the property
\begin{equation}\label{2.13}
L(x,y,Pz,Pw)=L(x,y,z,w),
\end{equation}
we call it a \emph{Riemannian $P$-tensor}.

If the curvature tensor $R$ on a Riemannian $\W_3$-manifold
$(M,P,g)$ is a Riemannian $P$-tensor, \ie{}
$R(x,y,Pz,Pw)=R(x,y,z,w)$, then $\tau^{**}=\tau$. Therefore
$\nP=0$, \ie $(M,P,g)$ is a Riemannian $P$-manifold.


\section{$P$-connection}

A linear connection $\n'$ on a Riemannian almost product manifold
$(M,P,g)$ preserving $P$ and $g$, \ie $\n'P=\n'g=0$, is called a
\emph{natural connection} \cite{Ga-Mi}.
\begin{defn}
The natural connection $\n'$ on a Riemannian almost product
manifold $(M,P,g)$ determined by
\begin{equation}\label{2.1}
    \n'_x y=\nabla_x y -\frac{1}{2}\bigl(\nabla_x P\bigr)Py,
\end{equation}
is called a \emph{$P$-connection}.
\end{defn}

Let $T$ be a torsion tensor of the $P$-connection $\n'$ determined
on $(M,P,g)$ by \eqref{2.1}. Because of the symmetry of $\nabla$,
from \eqref{2.1} we have $T(x,y)=-\frac{1}{2}\{\bigl(\nabla_x
P\bigr)Py- \bigl(\nabla_y P\bigr)Px\}$. Then, having in mind
\eqref{2.1}, we obtain
\[
    T(x,y,z)=g(T(x,y),z)=-\frac{1}{2}\{F(x,Py,z)-F(y,Px,z)\}.
\]
Hence and \eqref{3} we have
\begin{equation}\label{2.3}
    \mathop{\s} \limits_{x,y,z} T(x,y,Pz)=0.
\end{equation}

Let $Q$ be the tensor field determined by
\begin{equation}\label{2.4}
    Q(y,z)=-\frac{1}{2}\bigl(\nabla_y P\bigr)Pz.
\end{equation}
Having in mind \eqref{2}, for the corresponding (0,3)-tensor field
we have
\begin{equation}\label{2.5}
    Q(y,z,w)=-\frac{1}{2}F(y,Pz,w).
\end{equation}
Because of the properties \eqref{3}, \eqref{2.5} implies
$Q(y,z,w)=-Q(y,w,z)$.

Let $R'$ be the curvature tensor of the $P$-connection $\n'$.
Then, according to \eqref{2.1} and \eqref{5} we have \cite{Ko-No}
\begin{equation}\label{2.6}
\begin{array}{l}
    R'(x,y,z,w)=R(x,y,z,w)+\left(\nabla_x Q\right)(y,z,w)-\left(\nabla_y
    Q\right)(x,z,w)
  \\[4pt]
\phantom{R'(x,y,z,w)=}
       +Q\left(x,Q(y,z),w\right)-Q\left(y,Q(x,z),w\right).
\end{array}
\end{equation}

After a covariant differentiation of \eqref{2.5}, a substitution
in \eqref{2.6}, a use of \eqref{2}, \eqref{3}, \eqref{sigma} and
some calculations, from \eqref{2.6} we obtain
\[
\begin{array}{l}
    R'(x,y,z,w)=R(x,y,z,w)+\frac{1}{2}\left[\left(\nabla_x F\right)(y,z,Pw)-\left(\nabla_y
    F\right)(x,z,Pw)\right]
  \\[4pt]
\phantom{R'(x,y,z,w)=}
       +\frac{1}{4}\left[g\bigl(\left(\nabla_y P\right)z,\left(\nabla_x P\right)w\bigr)
       -g\bigl(\left(\nabla_x P\right)z,\left(\nabla_y P\right)w\bigr)\right].
\end{array}
\]

The last equality, having in mind the Ricci identity for
Riemannian almost product manifolds
\begin{equation*}
\bigl(\nabla_x F\bigr)(y,z,w)-\bigl(\nabla_y
F\bigr)(x,z,w)=R(x,y,Pz,w) - R(x,y,z,Pw),
\end{equation*}
implies
\begin{equation}\label{2.7}
\begin{array}{l}
    R'(x,y,z,w)\\[4pt]
    =\frac{1}{4}\left\{2R(x,y,z,w)+2R(x,y,Pz,Pw)+K(x,y,z,w)\right\},
\end{array}
\end{equation}
where $K$ is the tensor determined by
\begin{equation}\label{2.8}
    K(x,y,z,w)=-g\bigl(\left(\nabla_x P\right)z,\left(\nabla_y P\right)w\bigr)
       +g\bigl(\left(\nabla_y P\right)z,\left(\nabla_x
       P\right)w\bigr).
\end{equation}

In this way, the following theorem is valid.
\begin{thm}\label{thm2.2}
    The curvature tensor $R'$ of the $P$-connection $\n'$ on a
    Riemannian almost product manifold $(M,P,g)$ has the form \eqref{2.7}.\hfill$\Box$
\end{thm}

From \eqref{2.7} is follows immediately that the property
\eqref{2.11'} and \eqref{2.13} are valid for $R'$. Therefore, the
property \eqref{2.12} for $R'$ is a necessary and sufficient
condition $R'$ to be a Riemannian $P$-tensor. Since $R$ satisfies
\eqref{2.12}, then from \eqref{2.7} we obtain immediately the
following
\begin{thm}\label{thm2.3}
    The curvature tensor $R'$ of the $P$-connection $\n'$ on a
Riemannian $\W_3$-manifold $(M,P,g)$ is a Riemannian $P$-tensor
iff
\begin{equation}\label{2.9}
    2\mathop{\s} \limits_{x,y,z} R(x,y,Pz,Pw)=-\mathop{\s} \limits_{x,y,z}
    K(x,y,z,w).
\end{equation}\hfill$\Box$
\end{thm}

Let the following condition be valid for the Riemannian almost
product manifold $(M,P,g)$:

\begin{equation}\label{2.10}
    \mathop{\s} \limits_{x,y,z} R(x,y,Pz,Pw)=0.
\end{equation}

We say that the condition \eqref{2.10} characterizes a class
$\LL_2$ of the Riemannian almost product manifolds.

The equality \eqref{2.8} implies immediately the properties
\eqref{2.11'} and \eqref{2.13} for $P$. Then, according to
\eqref{2.9} and \eqref{2.10}, we obtain the following
\begin{thm}\label{thm2.4}
    Let $(M,P,g)$ belongs to the class $\LL_2$. Then the
    curvature tensor $R'$ of the $P$-connection $\n'$ is a Riemannian $P$-tensor iff
    the tensor $P$ determined by \eqref{2.8} is a Riemannian $P$-tensor, too.\hfill$\Box$
\end{thm}

Having in mind \eqref{2.7}, the last theorem implies the following
\begin{cor}\label{cor2.5}
    Let the curvature tensor $R'$ of the $P$-connection $\n'$ be
    a Riemannian $P$-tensor on $(M,P,g)\in\LL_2$. Then the tensor $H$,
    determined by
    \begin{equation}\label{2.11}
        H(x,y,z,w)=R(x,y,z,w)+R(x,y,Pz,Pw)
    \end{equation}
    is a Riemannian $P$-tensor, too.\hfill$\Box$
\end{cor}


\section{Curvature properties of the $P$-connection in $\W_3\cap\LL_2$}

Let us consider the manifold $(M,P,g)\in\W_3\cap\LL_2$ with a
Riemannian $P$-tensor of curvature $R'$ of the $P$-connection
$\n'$. Then, according to \thmref{thm2.4} and \corref{cor2.5}, the
tensors $K$ and $H$, determined by \eqref{2.8} and \eqref{2.11},
respectively, are also Riemannian $P$-tensors.

Let $\rho'$ and $\rho(K)$ be the Ricci tensors for $R'$ and $K$,
respectively. Then we obtain immediately from \eqref{2.7}
    \begin{equation}\label{3.1}
        \rho(y,z)+\rho^*(y,Pz)=2\rho'(y,z)-\frac{1}{2}\rho(K)(y,z).
    \end{equation}
From \eqref{3.1} we have
    \begin{equation}\label{3.2}
        \tau+\tau^{**}=2\tau'-\frac{1}{2}\tau(K),
    \end{equation}
where $\tau'$ and $\tau(K)$ are the scalar curvatures for $R'$ and
$K$, respectively.

It is known from \cite{Mek1}, that $\nP=2(\tau-\tau^{**})$. Then
\eqref{3.2} implies
    \begin{equation}\label{3.3}
        \tau=\tau'-\frac{1}{4}\left(\tau(K)-\nP\right).
    \end{equation}

From \eqref{2.8} we obtain
\[
    \rho(K)(y,z)=-g^{ij}g\bigl(\left(\nabla_{e_i} P\right)z,\left(\nabla_y
    P\right)e_j\bigr),
\]
from where
\[
    \tau(K)=g^{ij}g^{ks}g\bigl(\left(\nabla_{e_i} P\right)e_s,\left(\nabla_{e_k}
    P\right)e_j\bigr).
\]
Hence, applying \eqref{6}, we get
    \begin{equation}\label{3.4}
        \tau(K)=\frac{1}{2}\nP.
    \end{equation}
From \eqref{3.3} and \eqref{3.4} it follows
    \begin{equation}\label{3.5}
        \tau=\tau'+\frac{1}{8}\nP.
    \end{equation}

The equalities \eqref{3.2}, \eqref{3.3}, \eqref{3.4} and
\eqref{3.5} implies the following
\begin{prop}\label{prop3.1}
    Let the curvature tensor $R'$ of the $P$-connection $\n'$ be
    a Riemannian $P$-tensor on $(M,P,g)\in\W_3\cap\LL_2$.
    Then
    \begin{equation}\label{3.6'}
    \nP=-8(\tau'-\tau)=\frac{8}{3}(\tau'-\tau^{**})=2\tau(K).
    \end{equation}\hfill$\Box$
\end{prop}

\begin{cor}\label{cor4.2}
    Let the curvature tensor $R'$ of the $P$-connection $\n'$ be
    a Riemannian $P$-tensor on $(M,P,g)\in\W_3\cap\LL_2$.
    Then the following assertions are equivalent:
\begin{enumerate}
    \item[1)] $(M,P,g)$ a Riemannian $P$-manifold;
    \item[2)] $\tau'=\tau$;
    \item[3)] $\tau'=\tau^{**}$;
    \item[4)] $\tau(K)=0$.\hfill$\Box$
\end{enumerate}
\end{cor}

Let the considered manifold with a Riemannian $P$-tensor of
curvature $R'$ of the $P$-connection $\n'$ in $\W_3\cap\LL_2$ be
4-dimensional. Since $H$ is a Riemannian $P$-tensor, then
according to \cite{Sta}, we have
\begin{equation}\label{3.6}
      H=\nu(H)(\pi_1-\pi_2)+\nu^*(H)\pi_3,
\end{equation}
where $\nu(H)=\frac{\tau(H)}{8}$, $\nu^*(H)=\frac{\tau^*(H)}{8}$,
$\tau(H)$ and $\tau^*(H)$ are the scalar curvature of $H$ and its
associated one, and
\[
\begin{array}{l}
\pi_1(x,y,z,w)=g(y,z)g(x,w)-g(x,z)g(y,w),\\[4pt]
\pi_2(x,y,z,w)=g(y,Pz)g(x,Pw)-g(x,Pz)g(y,Pw),\\[4pt]
\pi_3(x,y,z,w)=g(y,z)g(x,Pw)-g(x,z)g(y,Pw),\\[4pt]
\phantom{\pi_2(x,y,z,w)=} +g(y,Pz)g(x,w)-g(x,Pz)g(y,w).
\end{array}
\]

From \eqref{2.7} and \eqref{2.11} we obtain
\begin{equation}\label{4.8}
\tau(H)=\frac{4\tau'-\tau(K)}{2},\qquad
\tau^*(H)=\frac{4\tau'^*-\tau^*(K)}{2},
\end{equation}
where $\tau'^*$ and $\tau^*(K)$ are the associated scalar
curvature to $\tau'$ and $\tau(K)$, respectively.

We apply \eqref{4.8} to \eqref{3.6} and thus we obtain the
following
\begin{prop}\label{prop3.2}
    Let the curvature tensor $R'$ of the $P$-connection $\n'$ be
    a Riemannian $P$-tensor on a 4-dimensional manifold $(M,P,g)\in\W_3\cap\LL_2$. Then
    \[
          H=\frac{4\tau'-\tau(K)}{16}(\pi_1-\pi_2)+\frac{4\tau'^*-\tau^*(K)}{16}\pi_3.
    \]\hfill$\Box$
\end{prop}

Let $\LL_1$ is the subclass of $\LL_2$ determined by
\begin{equation}\label{4.9}
    R(x,y,Pz,Pw)=R(x,y,z,w).
\end{equation}

The equalities \eqref{4.9} and \eqref{2.7} imply the following
\begin{prop}\label{prop3.3}
    Let $(M,P,g)\in\W_3\cap\LL_1$. Then
    \[
          R=R'-\frac{1}{4}K.
    \]\hfill$\Box$
\end{prop}

\begin{cor}\label{cor3.4}
    Let $(M,P,g)\in\W_3\cap\LL_1$. Then
    \[
         \tau=\tau'-\frac{1}{4}\tau(K),\qquad \tau^*=\tau'^*-\frac{1}{4}\tau^*(K).
    \]\hfill$\Box$
\end{cor}

\begin{cor}\label{cor3.5}
    Let $(M,P,g)\in\W_3\cap\LL_1$ and $\dim M=4$. Then
    \[
         \tau=\frac{1}{2}\tau(H),\qquad \tau^*=\frac{1}{2}\tau^*(H).
    \]\hfill$\Box$
\end{cor}



\bigskip

\textit{Dimitar Mekerov\\
Department of Geometry\\
Faculty of Mathematics and Informatics
\\
Paisii Hilendarski University of Plovdiv\\
236 Bulgaria Blvd.\\
4003 Plovdiv, Bulgaria
\\
e-mail: mircho@uni-plovdiv.bg}

\end{document}